\documentclass[a4paper,twoside]{article}

\def\shorttitle{Relative Chern Character and Supper-connection}
\def\shortauthor{Dexie Lin}

\usepackage{amsfonts}
\usepackage{amsmath}
\usepackage{amssymb}
\usepackage{amscd}
\usepackage{amsbsy}
\usepackage{amsthm}
\usepackage{bm}
\usepackage{empheq}
\usepackage{psfrag}
\usepackage{color}
\usepackage[colorlinks,linkcolor=blue,citecolor=red]{hyperref}
\usepackage{cite}
\usepackage{multicol}

\allowdisplaybreaks[4]

\newfont{\myfnt}{cmssi10 scaled 1440}
\numberwithin{equation}{section}
\catcode`@=11
\def\ps@nk{\def\@oddhead{\vbox{\hbox to \hsize{\pic \footnotesize \it \shorttitle
\hfill \rm \thepage} \vspace{1mm} \vspace*{-2mm}}}
\def\@evenhead{\vbox{\hbox to \hsize{\pic \footnotesize \rm \thepage \hfill \it \shortauthor}
\vspace{1mm} \vspace*{-2mm}}}
\def\@oddfoot{} \def\@evenfoot{}}
\def\ps@first{\def\@oddhead{\vbox{\hbox to \hsize{\pic \footnotesize
} \break}}
\def\@oddfoot{} \def\@evenfoot{}}
\newtheoremstyle{thmstyle}
  {6pt}
  {6pt}
  {\it}
  {}
  {\bf}
  {}
  {.5em}
  {}
\newtheoremstyle{remstyle}
  {6pt}
  {6pt}
  {\rm}
  {}
  {\bf}
  {}
  {.5em}
  {}
\def\Section#1{\Sec{\large #1} \setcounter{equation}{0} \vskip -6mm \indent}
\def\Sec{\@Startsection{section}{1}{\z@}
                                   {-3.5ex \@plus -1ex \@minus -.2ex}%
                                   {2.3ex \@plus.2ex}%
                                   {\normalfont\large\bfseries\boldmath}}
\def\@Startsection#1#2#3#4#5#6{%
  \if@noskipsec \leavevmode \fi
  \par
  \@tempskipa #4\relax
  \@afterindenttrue
  \ifdim \@tempskipa <\z@
    \@tempskipa -\@tempskipa \@afterindentfalse
  \fi
  \if@nobreak
    \everypar{}%
  \else
    \addpenalty\@secpenalty\addvspace\@tempskipa
  \fi
  \@ifstar
    {\@ssect{#3}{#4}{#5}{#6}}%
    {\@dblarg{\@Sect{#1}{#2}{#3}{#4}{#5}{#6}}}}
\def\@Sect#1#2#3#4#5#6[#7]#8{%
  \ifnum #2>\c@secnumdepth
    \let\@svsec\@empty
  \else
    \refstepcounter{#1}%
    \protected@edef\@svsec{\@seccntformat{#1}\relax}%
  \fi
  \@tempskipa #5\relax
  \ifdim \@tempskipa>\z@
    \begingroup
      #6{%
          \@hangfrom{\hskip #3\relax\@svsec \hskip -2.5mm}%
          \interlinepenalty \@M #8\@@par}
    \endgroup
    \csname #1mark\endcsname{#7}%
    \addcontentsline{toc}{#1}{%
      \ifnum #2>\c@secnumdepth \else
        \protect\numberline{\csname the#1\endcsname}%
      \fi
      #7}%
  \else
    \def\@svsechd{%
      #6{\hskip #3\relax
      \@svsec #8}%
      \csname #1mark\endcsname{#7}%
      \addcontentsline{toc}{#1}{%
        \ifnum #2>\c@secnumdepth \else
          \protect\numberline{\csname the#1\endcsname}%
        \fi
        #7}}%
  \fi
  \@xsect{#5}}
\renewenvironment{abstract}{%
        \small
        \quotation
         \noindent {\bfseries \abstractname } }%
      {\if@twocolumn\else\endquotation\fi}

\def\Subsec{\@StartSubsection{subsection}{2}{\z@}%
                                     {-3.25ex\@plus -1ex \@minus -.2ex}%
                                     {1.5ex \@plus .2ex}%
                                     {\normalfont\normalsize\bfseries\boldmath}}
\def\@StartSubsection#1#2#3#4#5#6{%
  \if@noskipsec \leavevmode \fi
  \par
  \@tempskipa #4\relax
  \@afterindenttrue
  \ifdim \@tempskipa <\z@
    \@tempskipa -\@tempskipa \@afterindentfalse
  \fi
  \if@nobreak
    \everypar{}%
  \else
    \addpenalty\@secpenalty\addvspace\@tempskipa
  \fi
  \@ifstar
    {\@ssect{#3}{#4}{#5}{#6}}%
    {\@dblarg{\@SubSect{#1}{#2}{#3}{#4}{#5}{#6}}}}
\def\@SubSect#1#2#3#4#5#6[#7]#8{%
  \ifnum #2>\c@secnumdepth
    \let\@svsec\@empty
  \else
    \refstepcounter{#1}%
    \protected@edef\@svsec{\@seccntformat{#1}\relax}%
  \fi
  \@tempskipa #5\relax
  \ifdim \@tempskipa>\z@
    \begingroup
      #6{%
          \@hangfrom{\hskip #3\relax\@svsec\hskip -1.5mm}%
          \interlinepenalty \@M #8\@@par}
    \endgroup
    \csname #1mark\endcsname{#7}%
    \addcontentsline{toc}{#1}{%
      \ifnum #2>\c@secnumdepth \else
        \protect\numberline{\csname the#1\endcsname}%
      \fi
      #7}%
  \else
    \def\@svsechd{%
      #6{\hskip #3\relax
      \@svsec #8}%
      \csname #1mark\endcsname{#7}%
      \addcontentsline{toc}{#1}{%
        \ifnum #2>\c@secnumdepth \else
          \protect\numberline{\csname the#1\endcsname}%
        \fi
        #7}}%
  \fi
  \@xsect{#5}}
\def\list#1#2{\ifnum \@listdepth >5\relax \@toodeep \else \global
\advance \@listdepth\@ne \fi \rightmargin \z@ \listparindent\z@
\itemindent\z@ \csname @list\romannumeral\the\@listdepth\endcsname
\def\@itemlabel{#1}\let\makelabel\@mklab \@nmbrlistfalse #2\relax
\@trivlist \parskip 0pt \parindent\listparindent \advance \linewidth
-\rightmargin \advance\linewidth -\leftmargin \advance\@totalleftmargin
\leftmargin \parshape \@ne \@totalleftmargin \linewidth \ignorespaces}
\renewcommand{\@makecaption}[2]{\begin{center}#1. #2\end{center}}
\catcode`@=12 
\theoremstyle{thmstyle}
\newtheorem{defi}{Definition}[section]
\newtheorem{thm}[defi]{\textbf{Theorem}}

\newtheorem{cor}[defi]{\textbf{Corollary}}
\newtheorem{lemma}[defi]{\textbf{Lemma}}

\newtheorem{prop}[defi]{Proposition}
\newtheorem{assumption}[defi]{Assumption}
\newsavebox{\mygraphic}
\def\pic{\begin{picture}(0,0) \put(-210,-1250){\usebox{\mygraphic}} \end{picture}}
\newfont{\HUGEbf}{cmbx10 scaled 3500}
\definecolor{gray}{rgb}{0.9,0.9,0.9}
\def\thebibliography#1{\section*{\bf \large References}
\list{[\arabic{enumi}]} {\settowidth \labelwidth{[#1]} \leftmargin
\labelwidth \advance \leftmargin \labelsep \usecounter{enumi}}
\def\newblock{\hskip .11em plus .33em minus .07em} \footnotesize \sloppy \clubpenalty
4000 \widowpenalty 4000 \sfcode`\.=1000 \relax}
\def\Ch{\mbox{Ch}}
\def\Tr{\mbox{Tr}}

\def\cs{\mbox{Cs}}

\def\spin{\mbox{spin}}

\def\ch{\mbox{ch}}

\def\cs{\mbox{cs}}

\theoremstyle{definition}

\numberwithin{equation}{section}

\title{Relative Chern character number and super--connection}
\author{Dexie Lin$^\dag$}

\date{}

\begin{document}

\maketitle

\thispagestyle{first}
\renewcommand{\thefootnote}{\fnsymbol{footnote}}

\footnotetext{\hspace*{-5mm} \begin{tabular}{@{}r@{}p{14cm}@{}} &
Manuscript last updated: \today.\\
$^\dag$ & Graduate School of Mathematical Sciences, The University of Tokyo, 3-8-1 Komaba, Meguro-ku, Tokyo 153-8914, Japan. E-mail: dexielin@ms.u-tokyo.ac.jp\\
\end{tabular}}

\renewcommand{\thefootnote}{\arabic{footnote}}

\begin{abstract}

For two complex vector bundles admitting a homomorphism, whose singularity locates in the disjoint union of some odd--dimensional spheres,
  we give a formula to compute the relative Chern characteristic number of these two complex vector bundles. In particular, for a $\spin$ manifold admitting some sphere bundle structure, we give a formula to express the index of a speccial twisted Dirac operator.

\vskip 4.5mm

\noindent\begin{tabular}{@{}l@{ }p{10cm}} {\bf Keywords } & Characteristic Classes , Super-connection, Chern character class\\
\end{tabular}

\vskip 4.5mm


\end{abstract}

\baselineskip 14pt

\setlength{\parindent}{1.5em}

\setcounter{section}{0}

\Section{Introduction}\label{sec-intro}

Let $M$ be a closed oriented  $2n$--dimensional  manifold.
The aim of this paper is to localize relative Chern character in the following setting. The relative Chern character is formulated as
\begin{equation}
  K(M,M\setminus U)\overset{\ch}{\longrightarrow} H^{even}(M,M\setminus U), \label{relative-Chern-character}
\end{equation}
where  $U$ is a proper open submanifold of $M$. In this paper, we suppose all manifolds are smooth.

We fix an element $[E]$ of $K(M,M\setminus U)$ which is represented by the following data  \cite[Lemma 8.4]{ABS}:
$$0\to E_+\overset{v}{\longrightarrow}E_-\to0,$$
where $E_+$ and $E_-$ are two complex vector bundles with the same rank, and $v$ is a complex
linear homomorphism between the two vector bundles which is isomorphic on $M\setminus U$.
We call $x\in M$ a singularity point of $v$, if $v_x$ is not isomorphic. The singularity of $v$ is the set of singular points, denoted by $Sing(v)$.
We assume that $U$ is a tubular neighborhood of a closed submanifold $Y$.
 If we deform $v$ slightly,
we can assume that $Sing(v)$ is a subset of $Y$. In this paper we
assume $Sing(v) \subset Y$. Note that we do not assume the
transversality of $v$.

 The pairing  between the  relative Chern character class $\ch(E_+,E_-,v)$ and the fundamental class $[M]$ is equal to
$\langle\ch(E_+)-\ch(E_-),[M]\rangle$, which  we call  relative Chern character number.
The relative Chern character number depends only on the data $(E_+,E_-,v)$ restricted on the tubular neighborhood of  $Y$. However, the explicit relationship between the relative  Chern character number and the neighborhood is not clear, in particular when $v$ is not transverse to zero as a section of $\hom(E_+,E_-)$.

The notion of super--connection was introduced by Quillen \cite{Q}, as a generalization of the notion of connection in the category of $\mathbb Z_2$--graded vector bundles.
We can use the super--connection to localize the relative Chern character number. The following examples are typical  localizations of the relative Chern character in different conditions on $Y$ by super--connection.

\begin{itemize}
   \item[1)] When $Y$ is a set consisting of finitely many points,  Feng, Li and Zhang \cite{FLZ} gave a   Poincar\'e-Hopf type formula to calculate  the relative Chern character number.
   \item[2)] When $M$ and $Y$ are  closed complex manifolds, $Y\hookrightarrow M$ is a holomorphic embedding and $Sing(v)=Y$, Bismut(1990) gave several estimates of the  convergence of the relative Chern character with respect to the deformation variable of the super--connection, under the "quantization assumption"
   \cite[Section 1 a. Assumptions and notaions]{B90}. One of the estimates implies the localization of the relative Chern character number, and the number turns out to be equal to some topological invariants on $Y$.
 \end{itemize}
We suppose that $v$ satisfies the following  assumption.
\begin{assumption}\label{Basic-Assumption-of-paper}
There are finite embeddings $\iota_i:S^{2k_i-1}\hookrightarrow M$, $i=1,\cdots,l$ for some integers $1\leq k_1,\cdots,k_l< n$, such that
\begin{itemize}
  \item for each $S^{2k_i-1}$, its normal bundle is trivial,
  \item the images of these embeddings are mutually disjoint, i.e. $\iota_i(S^{2k_i-1})\cap
  \iota_j(S^{2k_j-1})=\emptyset$, for $i\neq j$.
\end{itemize}
  We suppose that  $Sing(v)$ is a subset of the union of these images, i.e. $Sing(v)\subset\bigcup_{1\leq i\leq l} Im(\iota_i)$.
\end{assumption}
For  convenience, we also  use $S^{2k_i-1}$ to denote  $Im(\iota_i)$. Let $N(S^{2k_i-1})$ denote the tubular neighborhood of $S^{2k_i-1}$. By the above assumption, we have a diffeomorphism $N(S^{2k_i-1})\cong S^{2k_i-1}\times D^{2n-2k_i+1}$, where $D^{2n-2k_i+1}$ denotes the standard unit disk of dimension $2n-2k_i+1$, and the boundary $\partial N(S^{2k_i-1})$ is diffeomorphic to $S^{2k_i-1}\times S^{2n-2k_i}$.

When $Sing(v)$ is just a set of finitely many points, Feng, Li and Zhang \cite{FLZ} gave an explicit formula of  the relative Chern character number as a sum of  contributions from each point $p\in Sing(v)$, which depends only on the restriction of $v$ on the boundary of a small disk $D(p)$ around $p$. Actually the homomorphism
 $$v\big|_{\partial D(p)}:E_+\big|_{\partial D(p)}\to E_-\big|_{\partial D(p)}$$
  determines an element in $K^1(S^{2n-1})\cong\mathbb Z$. They used the notion of  odd Chern character due to  Getzler \cite[Section 1,P. 492]{G} to calculate the contribution.

 In this paper, we follow the strategy of Feng-Li-Zhang's work combined with  a calculation of some homotopy set which will be stated in Lemma  \ref{lemma-split.of.homotopy.group}. Our method to localize the relative Chern character number is to use the odd Chern character on the product sphere, i.e. $\partial N(S^{2k_i-1}) \cong S^{2k_i-1} \times S^{2n-2k_i}$. In the rest of  this paper, we fix these diffeomorphisms. We will show  that
the homomorphism
$$
v\big|_{\partial N(S^{2k_i-1})}:E_+\big|_{\partial N(S^{2k_i-1})}\to E_-\big|_{\partial N(S^{2k_i-1})}$$
also determines an integer, which we will denote by $\deg^*(v_{i})$.

 The main theorem of this paper is:

\begin{thm}\label{Main-Theorem}
  Suppose $M$ is a closed oriented $2n$--dimensional  manifold, then
  under the assumption \ref{Basic-Assumption-of-paper}, the following formula holds.
  \begin{equation}
  \langle \ch(E_+)-\ch(E_-),[M]\rangle=(-1)^{n+1}\sum_{1\leq i\leq l}\deg^*(v_{i})\label{Main-Degree-Formula}
  \end{equation}
\end{thm}

 A direct application of the above formula is the localization of the index of twisted Dirac operator for some $\spin$ manifold.
 To be more specific, let $M$ be a closed oriented  $2n$--dimensional $\spin$ manifold, $\pi:M\to B$ be a bundle with sphere fiber of dimension $2k-1$, for some integer $2\leq k<n$  and $F$ be a complex vector bundle of rank $r$ over $M$.
We consider the case that $E_+$ is a trivial rank $r$ bundle and
$E_-$ is equal to $F$. Assume that $v$ is a homomorphism between $E_+$  and $E_-$ satisfying $Sing(v)\subset\coprod_{1\leq i\leq l}\pi^{-1}(x_i)$, where $x_1,\cdots,x_l$ are finitely many points in $B$, such that the  homology class $[\pi^{-1}(x_i)]$ belongs to the torsion part of $H_{2k-1}(M;\mathbb Z)$ for each $x_i$. Let $D^{\otimes F}$ be the twisted Dirac operator, i.e. Dirac operator of the $F$--tensered spinor bundle. Then, the index can be localized as follows(see Proposition \ref{proposition-index-localization}):
 \begin{equation*}
 Ind(D^{\otimes F}_+)=(-1)^n\sum_{1\leq i\leq l}
 \deg^*(v_{x_i}).\label{Index-Localization-Formula}
 \end{equation*}

The organization of the paper is as follows:

 In  Section \ref{sec-pre}, we  review some properties of odd Chern characteristic class. In Section \ref{sec-local.calculus}, we give the local calculation near the singularities contained in the product of two spheres. In Section \ref{sec-proof}, we give the proof of our Main  Theorem. In Section \ref{sec-application}, we show how to apply such a formula to localize the index of the twisted Dirac operator for the spin manifold mentioned above.

\Section{Review of odd Chern character}\label{sec-pre}

Let $X$ be a closed manifold.
We recall the notion of odd Chern character given by Getzler \cite[P.490-495]{G}.

If $\nabla_0$ and $\nabla_1$ are two connections on a complex vector bundle $E$ over $X$, their Chern-Simons form is the differential form defined by
\[\cs(\nabla_0,\nabla_1)=\int^1_0\Tr(\dot\nabla_ue^{\nabla^2_u})du,\]
where $\nabla_u=(1-u)\nabla_0+u\nabla_1$.
By the projection
$[0,1]\times X\to X$, we can define the connection on the product space $[0,1]\times X$, $$\tilde\nabla=du\frac\partial{\partial u}+\nabla_u.$$
Let $N$ be a larger enough integer.

\begin{defi}
  The odd Chern character $\Ch(g)$ of a differentiable map $g:X\to GL_N(\mathbb C)$ is $\cs(d,d+g^{-1}dg)$.
\end{defi}

\begin{prop}[cf. {\cite[Proposition 1.2]{G}}]
  The odd Chern character is a closed differential form, and can be formulated as the following:
  \[\Ch(g)=\sum^\infty_{k=0}(-1)^k\frac{k!}{(2k+1)!}\Tr(\omega^{2k+1}),\]
  where $\omega=g^{-1}dg$.
\end{prop}

Let $g_t$, $t\in[0,1]$, be a family of differential maps from $X$ to $GL_N(\mathbb C)$. It defines a differential mapping $\tilde g:[0,1]\times X\to  GL_N(\mathbb C)$. Hence, the odd Chern character $\Ch(\tilde g)$ can be decomposed into the following:
\[\Ch(\tilde g)=\Ch(g_t)+dt\wedge\widetilde\Ch(g_t).\]

\begin{prop}[cf. {\cite[Proposition 1.3]{G}}]
  $\widetilde\Ch(g_t)$ can be expressed by the formula
  \[\widetilde\Ch(g_t)=\sum^\infty_{k=0}(-1)^k\frac{k!}{(2k)!}\Tr(g^{-1}_t\dot g_t\wedge
  \omega^{2k}_t),\]
  and satisfies the transgression formula
  \[\frac\partial{\partial t}\Ch(g_t)=d\widetilde\Ch(g_t).\]
\end{prop}

The above  proposition implies that the cohomology class $[\Ch(g)]$ depends only on the homotopy class $[g]\in[X,GL_N(\mathbb C)]$.

\begin{lemma}[cf. {\cite[Proposition 1.4]{G}}]\label{lemma-integral.on.sphere}
  If $g:S^{2k-1}\to GL_N(\mathbb C)$ is a differential mapping, then we have the integral
  \begin{equation}
    \frac1{(-2\pi\sqrt{-1})^{k}}\int_{S^{2k-1}}\Ch(g)=-\deg^{top}(g),\label{degree2}
  \end{equation}
  where $\deg^{top}(g)$ denotes the topological degree of the homotopy class $[g]\in\pi_{2k-1}(GL_N(\mathbb C))$.
\end{lemma}

In this paper, we denote the integral $\frac1{(-2\pi\sqrt{-1})^{k}}\int_{S^{2k-1}}\Ch(g)$ by $\deg(g)$.

\Section{Local calculation on singularity}\label{sec-local.calculus}

Now we consider the local model of each $S^{2k_i-1}$. To simplify the notation, we  drop the subscript $i$.
 We need the following  technique in the homotopy theory to split the elements in $[S^{2n-2k}\times S^{2k-1},GL_N(\mathbb C)]$.

\begin{lemma}\label{lemma-split.of.homotopy.group}
  Any homotopy class $[g]\in[S^{2n-2k}\times S^{2k-1}, GL_N(\mathbb C)]$ has a representation as a matrix product of two mappings $f:S^{2k-1}\to GL_N(\mathbb C)$ and $h:S^{2n-1}\to GL_N(\mathbb C)$, i.e. there exists $g\in[g]$, such that $g=pr^*_2f\cdot \phi^*h$, where $pr_2$ denotes the canonical projection $S^{2n-2k}\times S^{2k-1}\to S^{2k-1}$ and $\phi$ denotes the smash product $S^{2n-2k}\times S^{2k-1}\overset{\wedge}{\to}S^{2n-1}$.
\end{lemma}

\begin{pf}
  For any positive integers $p$ and $q$, we have the cofibration sequence,
  $$S^{p+q-1}\to S^p\vee S^q\to S^p\times S^q\to S^{p+q}\to S^{p+1}\vee S^{q+1},$$
  where the first map denotes the attaching and the forth map denotes the suspension, therefore null-homotopic.
  Taking $p=2n-2k$ and $q=2k-1$, we have the exact short sequence
  $$1\to\pi_{2n-1}(GL_N(\mathbb C))\to
  [S^{2n-2k}\times S^{2k-1},GL_N(\mathbb C)]\to
  \pi_{2n-2k}(GL_N(\mathbb C))\oplus\pi_{2k-1}(GL_N(\mathbb C))\to1.$$
  Since $\pi_{2n-2k}(GL_N(\mathbb C))=0$, we can find a homotopic equivalent map $$g_1:S^{2n-2k}\times S^{2k-1}\to GL_N(\mathbb C),\mbox{ such that }
  g_1\Big|_{S^{2n-2k}\times q}=Id_N,$$
  where $q$ is a fixed point in $S^{2k-1}$.
  Fixing a point $p$ in $S^{2n-2k}$, we set $f=g_1\Big|_{p\times S^{2k-1}}$.
  By the projection $p_2:S^{2n-2k}\times S^{2k-1}\to \{p\}\times S^{2k-1}$, we have that $(p^*_2f)^{-1}g_1$ is null-homotopy in the above exact sequence, hence it equals to a map  $h:S^{2n-1}\to GL_N(\mathbb C)$.
\end{pf}

By Lemma \ref{lemma-integral.on.sphere} and Lemma \ref{lemma-split.of.homotopy.group}, we can calculate the integral on the product sphere.

\begin{prop}
  If $[g]\in[S^{2n-2k}\times S^{2k-1},GL_N(\mathbb C)]$, then
  \begin{equation}
    \frac1{(-2\pi\sqrt{-1})^n}\int_{S^{2n-2k}\times S^{2k-1}}\Ch(g)=\deg(h), \label{product-degree}
  \end{equation}
  where $h$ is  a differential mapping from  $S^{2n-1}$ to $GL_N(\mathbb C)$ as in the above lemma.
\end{prop}

Here, we denote the integral $\frac1{(-2\pi\sqrt{-1})^n}\int_{S^{2n-2k}\times S^{2k-1}}\Ch(g)$ by $\deg^*(g)$.

\begin{pf}
    By Lemma \ref{lemma-split.of.homotopy.group}, we can find two differential mappings $f:S^{2k-1}\to GL_N(\mathbb C)$ and $h:S^{2n-1}\to GL_N(\mathbb C)$, such that there exists a representation $g:S^{2n-2k}\times S^{2k-1}\to GL_N(\mathbb C)$ of $[g]$ as the following
   $$g=pr^*_2f\cdot\phi^*h,$$
   where
$\phi$ denotes the smash product
  $$S^{2n-2k}\times S^{2k-1}\overset{\wedge}{\to}S^{2n-2k}\wedge S^{2k-1}=S^{2n-1}$$
     and $pr_2$ denotes the canonical projection
    $$S^{2n-2k}\times S^{2k-1}\to S^{2k-1}.$$
  It suffices to show that the integral is independent on  $f$.
  Note that we have a canonical embedding $GL_N(\mathbb C)\hookrightarrow GL_{2N}(\mathbb C)$, the image of $pr^*_2f\cdot \phi^*h$ under this embedding is
   $$\left(\begin{array}{cc}
   pr^*_2f\cdot \phi^*h&\\
   &Id_N\\
   \end{array}\right)\in GL_{2N}(\mathbb C),$$
   which is homotopic to $$\left(\begin{array}{cc}
   pr^*_2f&\\
   &\phi^*h\\
   \end{array}\right).$$
   Hence, we have $$[\Ch(g)]=[\Ch(pr^*_2f\cdot \phi^*h)]=[pr^*_2\Ch(f)]+[\phi^*\Ch(h)].$$
   Moreover, we know that $2k-1<2n-1$ and the differential form $\Ch(f)$  only depends on the  component $S^{2k-1}$, therefore the integral $\int_{S^{2n-2k}\times S^{2k-1}}pr^*_2\Ch(f)$ is zero.

   Now, we obtain
   \begin{eqnarray*}
     \frac1{(-2\pi\sqrt{-1})^n}\int_{S^{2n-2k}\times S^{2k-1}}\Ch(g)&=&
    \frac1{(-2\pi\sqrt{-1})^n}\int_{S^{2n-2k}\times S^{2k-1}}\phi^*\Ch(h)+pr^*_2\Ch(f)\\
    &=& \frac1{(-2\pi\sqrt{-1})^n}\int_{S^{2n-2k}\times S^{2k-1}}\phi^*\Ch(h)\\
    &=& \frac1{(-2\pi\sqrt{-1})^n}\int_{\phi(S^{2n-2k}\times S^{2k-1})}\Ch(h)\\
    &=&  \frac1{(-2\pi\sqrt{-1})^n}\int_{S^{2n-1}}\Ch(h)=\deg(h).
   \end{eqnarray*}
\end{pf}

\Section{Proof of  Theorem \ref{Main-Theorem}}\label{sec-proof}

Without loss of generality, we suppose that the rank of  $E_\pm$ is  large enough.

Let $E=E_+\oplus E_-$ be the $\mathbb Z_2$--graded complex vector bundle over $M$. Choose two unitary connections $\nabla^{E_+}$ and $\nabla^{E_-}$ on $E_+$ and $E_-$ respectively, such that near each small tubular neighborhood of  $S^{2k_i-1}$, the connections are both trivial. Let $\nabla^E=\nabla^{E_+}\oplus\nabla^{E_-}$ be the $\mathbb Z_2$--graded connection on $E$, $v:E_+\to E_-$ extend to an endomorphism of $E$ by acting as zero on $E_-$ and $v^*$ be the adjoint of $v$ with respect to the hermitian metrics on $E_\pm$ respectively. Set $V=v+v^*$, then $V$ is an odd endomorphism of $E$ and $V^2$ is fiberwise positive on $M_1=M\setminus \bigcup\limits_i N(S^{2k_i-1})$.

For the convenience of computation, we define a function $\varphi: \Omega^*(M)\to \Omega^*(M)$, by
\[\varphi(\alpha)=(2\pi\sqrt{-1})^{-\frac k2}\alpha,\]
where $\alpha$ is a differential form of degree $k$.

For any $T\in\mathbb R^{\geq0}$, let $A_T$ be the super--connection on $E$ in the sense of Quillen(cf. \cite[Section 2]{Q}), defined by
\begin{equation}
  A_T=\nabla^E+TV.
\end{equation}
Let $\ch(E,A_T)$ be the associated Chern character form defined by
\begin{equation}
  \ch(E,A_T)=\varphi(\Tr[e^{-A^2_T}]).
\end{equation}
By  straightforward calculation, one can derive the following lemmas.

\begin{lemma}[cf. {\cite[(2.6) and (2.5)]{FLZ}}, {\cite[Proposition 2]{Q}}]
  \begin{equation}
    \frac{\partial\ch(E,A_T)}{\partial T}=-\frac1{\sqrt{2\pi\sqrt{-1}}}d\varphi
    \Tr[Ve^{-A^2_T}].\label{transgression}
  \end{equation}
\end{lemma}
\noindent We set
\begin{equation}
  \gamma(T)=\frac1{\sqrt{2\pi\sqrt{-1}}}\varphi\int^T_0\Tr(Ve^{-A^2_t})dt.
\end{equation}
Then, we get the transgression  formula

\begin{equation}
  \ch(E,A_0)-\ch(E,A_T)=d\gamma(T).\label{transgression3}
\end{equation}

\begin{lemma}[cf. {\cite[Section 4]{Q}}]
 On $M_1$, we have
  \begin{equation}
    \lim_{T\to+\infty}\ch(E,A_T)=0.\label{vanishing}
  \end{equation}
\end{lemma}

\begin{lemma}[cf. {\cite[Lemma 2.1]{FLZ}}]
  The relative Chern character number can be expressed as follows:
  \begin{equation}
    \langle\ch(E_+)-\ch(E_-),[M]\rangle=-\sum_{1\leq i\leq l}\lim_{T\to\infty}\int_{\partial N(S_{x_i})}\gamma(T).\label{bound}
  \end{equation}
\end{lemma}

Now, we can give the proof of Theorem \ref{Main-Theorem}.

\begin{pf}
  When $v$ is restricted on  $\partial N(S^{2k_i-1})$,  we get that $v^*=v^{-1}$ and $V^2$ is the identity map acting on $E\mid_{\partial N(S^{2k_i-1})}$. Since $\nabla^E$ is  trivial over $N(S^{2k_i-1})$ which is denoted  by $d$, then on $\partial N(S^{2k_i-1})$ we have the following  two identities
  \[A_t=d+tV,\quad A^2_t=d^2+t^2V^2+t[d,V]=t^2Id_E+tdV.\]
  Then,   we get
  \begin{eqnarray*}
    \int_{\partial N(S^{2k_i-1})}\gamma(T)&=&\frac1{\sqrt{2\pi\sqrt{-1}}}
    \int_{\partial N(S^{2k_i-1})}\varphi \int^T_0\Tr(Ve^{-A^2_t})dt\\
    &=&\frac1{\sqrt{2\pi\sqrt{-1}}}\int_{\partial N(S^{2k_i-1})}\varphi\int^T_0e^{-t^2}
    \Tr(Ve^{-tdV})dt
  \end{eqnarray*}
  Since $e^{-tdV}=e^{-t\left(\begin{array}{cc}
  &dv^*\\
  dv&\\
  \end{array}\right)}$, after taking the super--trace only the odd degree terms of $\sum\limits_n\left(\begin{array}{cc}
  &dv^*\\
  dv&\\
  \end{array}\right)^n$ remain. By the identity
  \[\left(\begin{array}{cc}
  &dv^*\\
  dv&\\
  \end{array}\right)^{2k-1}=\left(\begin{array}{cc}
  &dv^*(dvdv^*)^{k-1}\\
  dv(dv^*dv)^{k-1}&\\
  \end{array}\right),\]
  we have
  \begin{eqnarray*}
    &&\int_{\partial N(S^{2k_i-1})}\varphi\int^T_0e^{-t^2}\Tr(Ve^{-tdV})dt\\
    &=&\sum_{k\geq1}\int^T_0e^{-t^2}t^{2k-1}\frac{-1}{(2k-1)!}
    \int_{\partial N(S^{2k_i-1})}\varphi(\Tr_{E_+}(v^*dv(dv^*dv)^{k-1})-\Tr_{E_-}(vdv^*(dvdv^*)^{k-1}))\\
    &=&\frac1{(2\pi\sqrt{-1})^{n-\frac12}}\frac{2(-1)^n}{(2n-1)!}\int^T_0t^{2n-1}e^{-t^2}dt
    \int_{\partial N(S^{2k_i-1})}\Tr_{E_+}(v^{-1}dv)^{2n-1}.
  \end{eqnarray*}

  Therefore, we get
  \begin{eqnarray*}
    \lim_{T\to+\infty}\int_{\partial N(S^{2k_i-1})}\gamma(T)&=&
    \frac1{(2\pi\sqrt{-1})^n}\frac{2(-1)^n}{(2n-1)!}\int^{+\infty}_0t^{2n-1}e^{-t^2}dt
    \int_{\partial N(S^{2k_i-1})}\Tr_{E_+}(v^{-1}dv)^{2n-1}\\
    &=&\frac1{(2\pi\sqrt{-1})^n}\frac{(-1)^n(n-1)!}{(2n-1)!}
    \int_{\partial N(S^{2k_i-1})}\Tr_{E_+}(v^{-1}dv)^{2n-1}\\
    &=&(-1)^{n}\deg^*(v_i).
  \end{eqnarray*}
  The last equality follows from the formula \eqref{product-degree}.
  After taking the sum of $i\in\{1,\cdots, l\}$, we obtain the  formula \eqref{Main-Degree-Formula}.
\end{pf}

At the end of this section, we want to point out that the above argument can also be applied to the case, when $Sing(v)$ is a set consisting of finite points. In other words,   Theorem \ref{Main-Theorem} is a generalization of Feng-Li-Zhang's result.

\begin{thm}
  [Feng, Li and Zhang {\cite[Theorem $1$]{FLZ}} ]
  Let $M$ be a closed oriented  manifold of dimension $2n$, and $E_+$, $E_-$ be two complex vector bundles with same rank. Given a homomorphism $v\in\Gamma(\hom(E_+,E_-))$, whose singularity set consists of finite points of $M$, then
  $$\langle\ch(E_+)-\ch(E_-),[M]\rangle=(-1)^{n-1}\sum_{p\in Sing(v)}\deg(v_p),$$
  where $\deg(v_p)$ denotes the degree of the map $v\big|_{\partial D(p)}:E_+\big|_{\partial D(p)}\to E_-\big|_{\partial D(p)}$.
\end{thm}

\Section{An application to the index of twisted Dirac operator}\label{sec-application}

In this section, we  show that  the index of some twisted Dirac operator can be localized to the relative Chern character number. We  assume that $M$ is a closed oriented $2n$--dimensional spin manifold, and satisfies the following assumption.

\begin{assumption}\label{Basic-Assumption-of-Manifold}   $M$ admits the odd--dimensional sphere bundle structure over a closed manifold $B$, i.e. there is a smooth map
$$\pi:M\overset{S}\longrightarrow B,$$
where $B$ and $S$ are closed manifolds, such that  the following conditions hold: \begin{itemize}
\item[1)]$M$ is locally split, i.e.  each $x\in B$ has a small neighborhood $U(x)$ in $B$, such that $\pi^{-1}(U(x))$ is diffeomorphic to $U(x)\times S$,
\item[2)]$S$ is diffeomorphic to an odd dimensional sphere $S^{2k-1}$, for some integer $2\leq k<n$.\end{itemize}
\end{assumption}

For technical reason, we need the following lemma.

\begin{lemma}\label{lemma-local-exact}
  Let $X$ and $Y$ be two closed manifolds, and $\iota$ denote the embedding $Y\hookrightarrow X$.
  Suppose $Y$ is a rational homology sphere and $[Y]\in Tor(H_*(X;\mathbb Z))$.
  Then,  we have that for any $d$--closed form $\omega$ on $X$ without the zero degree part, the restriction  $\omega\mid_{\mathcal N}$ is $d$--exact, where $\mathcal N$ denotes  a small tubular neighborhood of $Y$.
\end{lemma}

\begin{pf}
  We can choose a tubular neighborhood of $Y$ by equipping a metric on $X$, and by the pull-back of the projection
  $\mathcal N\to Y$, we need only to show that $\omega\mid_{Y}$ is $d$-exact.
  By the Poincar\'e duality, it  suffices  to show  that the following formula holds
  \begin{eqnarray*}
    0=\langle\omega,\iota_* [Z]\rangle,
  \end{eqnarray*}
  for any element  $[Z]\in H_*(Y;\mathbb Z)$.
  Since $Y$ is a rational homology sphere and $[Y]\in Tor(H_*(X;\mathbb Z))$, the above formula holds.
\end{pf}

\begin{prop}\label{proposition-index-localization}
  Let $M$ be a closed  oriented $2n$--dimensional $\spin$ manifold  satisfying the assumption \ref{Basic-Assumption-of-Manifold}.  Suppose there is a linear homomorphism
  $\underline{\mathbb C}^r\overset{v}{\to}F$ between the two complex vector bundles $\underline{\mathbb C}^r$ and $F$ with the same rank $r$, satisfying that $Sing(v)\subset\coprod_{1\leq i\leq l}
  \pi^{-1}(x_i)$ for finitely many points $\{x_i\}_{1\leq i\leq l}$ in $B$ and $[\pi^{-1}(x_i)]\in Tor(H_{2k-1}(M;\mathbb Z))$ for each $x_i$, then
  $$Ind(D^{\otimes F}_+)=(-1)^n\sum_{1\leq i\leq l}\deg^*(v_{x_i}).$$
\end{prop}

\begin{pf}
  Since $S$ can be viewed as the foliation of the manifold $M$, and it admits a positive scalar curvature. By Connes Vanishing Theorem \cite[Theorem $02$]{Connes} we know $\langle\hat{A}(TM),[M]\rangle=0$. We omit the proof here and refer to a geometric proof \cite{Zh15}, which is given by   Zhang.
  Actually, in our setting we can use adiabatic limits method of Dirac operator as in Bismut's work \cite[Proposition $5.2$]{B86}  to show that $\hat{A}(TM)$ vanishes.

  Let $E:=\underline{\mathbb C}^r\oplus F$ be the $\mathbb Z_2$--graded vector bundle.
   We need to show that
$$\int_M\hat{A}(TM)(\ch(E))=\int_M\ch(E).$$
Since we can write
$$\hat{A}(TM)=\sqrt{\det\left(\frac{R^{TM}/2}{\sinh(R^{TM}/2)}\right)}=1+\beta,$$
where $R^{TM}$ denotes the curvature of $TM$ and $\beta$ denotes the  form with high degree,
it suffices to  show a more general formula:
$$\langle\omega\wedge\ch(E),[M]\rangle=\langle\omega^0\wedge\ch(E),[M]\rangle,\quad \mbox{ for any }d-\mbox{closed form } \omega,$$
where $\omega^0$ denotes the $0$-degree part of $\omega$.

Let $N_\epsilon(\pi^{-1}(x_i))$ be the tubular neighborhood of $\pi^{-1}(x_i)$ of radius $\epsilon$ in normal direction and $\mathcal U_\epsilon=\coprod_{1\leq i\leq l}
  N_\epsilon(\pi^{-1}(x_i))$. Find a cut-off function whose support locates in $\mathcal U_{2\epsilon}$ and satisfying $\rho\mid_{\mathcal U_\epsilon}\equiv 1$.
By Lemma \ref{lemma-local-exact},   there is a form $\alpha$ on $\mathcal U_{2\epsilon}$  such that  $\omega^*= d\alpha$ in $\mathcal U_\epsilon$, where $\omega^*$ denotes $\omega-\omega^0$.
 Hence, we get
\begin{eqnarray*}
&&\langle\omega^*\wedge\ch(E),[M]\rangle-\langle d (\rho\cdot\alpha)\wedge\ch(E),[M]\rangle\\
&=&\int_M\left((1-\rho)\omega^*-d\rho\wedge\alpha\right)\wedge\ch(E)+
\int_M\rho\omega^*\wedge\ch(E)-\rho d\alpha\wedge\ch(E)\\
&=&\int_{M\setminus \mathcal U_\epsilon}(\omega^*-d(\rho\alpha))\wedge\ch(E,A_T)+\int_{\mathcal U_\epsilon}(\rho\omega^*-d\alpha)\wedge\ch(E).
\end{eqnarray*}
 By the formula \eqref{vanishing} on $M\setminus \mathcal U_\epsilon$ and $\rho\omega^*=d\alpha$ on $\mathcal U_\epsilon $, we  say that the above formula equals zero.
 Therefore, we get
 \begin{eqnarray*}
  Ind(D^{\otimes F}_+)&=&-\left(\int_M\hat{A}(TM)\ch(\underline{\mathbb C}^r)-\int_M\hat{A}(TM)\ch(F)\right)\\
  &=&-\int_M\ch(E)=(-1)^n\sum_{1\leq i\leq l}\deg^*(v_{x_i}).
 \end{eqnarray*}
\end{pf}

\begin{cor}
  Let $M$ be a closed oriented $2n$--dimensional manifold satisfying the assumption \ref{Basic-Assumption-of-Manifold}. If a complex line bundle $L$ over $M$ admits a global section $s\in\Gamma(M,L)$, such that there are finitely many points $\{x_i\}_{1\leq i\leq l}$ in $B$ satisfying  $s^{-1}(0)\subset\coprod_{1\leq i\leq l}\pi^{-1}(x_i)$  and $[\pi^{-1}(x_i)]\in Tor(H_{2k-1}(M;\mathbb Z))$ for each $x_i$, then
  \begin{equation}
  Ind(D^{\otimes L}_+)=(-1)^n\sum_{1\leq i\leq l}\deg^*(s_{x_i}).\label{IndexFormula}
  \end{equation}
\end{cor}

\begin{pf}
 The global section $s$ can be viewed as the complex  homomorphism between the trivial line bundle and $L$, i.e. we have a homomorphism
  $$\underline{\mathbb C}\overset{s}{\to}L.$$
 By the similar arguments  in the  proof  of Proposition \ref{proposition-index-localization}, we have
  $$\begin{array}{c}
    -Ind(D^{\otimes L}_+)=\int_M\hat{A}(TM)(\ch(\underline{\mathbb C})-\ch(L))=(-1)^{n+1}\sum\limits_{1\leq i\leq l}
    \deg^*(s_{x_i})\\
    ~Ind(D^{\otimes L}_+)=(-1)^n\sum\limits_{1\leq i\leq l}
    \deg^*(s_{x_i}).\\
  \end{array}$$
\end{pf}

Up to here we have considered the cases for odd-dimensional singularity. However, we also have the next corollary for a special case in even--dimensional singularity.

\begin{cor}
  Let $\pi:M\to B$ be a bundle with fiber $S^{2n-2}$, where $M$ is a closed oriented spin manifolds of dimension $2n$ for $n\geq3$ and $B$ is a closed manifold of dimension  $2$. Suppose that $L$ is a complex line bundle over $M$ and $s\in\Gamma(M,L)$, such that there are finitely many points $\{x_i\}_{1\leq i\leq l}$ in $B$ satisfying  $s^{-1}(0)\subset\coprod_{1\leq i\leq l}\pi^{-1}(x_i)$ and the fundamental class $[\pi^{-1}(x_i)]$ belongs to the torsion part of $H_{2n-2}(M;\mathbb Z)$, for each $x_i$. Then, we have
  $$Ind(D^{\otimes L}_+)=(-1)^n\sum_{1\leq i\leq l}\deg^*(s_{x_i}).$$
\end{cor}

\begin{pf}
   We  know that the normal bundle over $\pi^{-1}(x_i)$ is trivial and we can choose a connection of $L$, such that near $\pi^{-1}(x_i)$ it is trivial, because any complex line bundle over $\pi^{-1}(x_i)\cong S^{2n-2}$ is trivial for $n\geq3$.
   The tubular neighborhood of $\pi^{-1}(x_i)$ is diffeomorphic to $D^2 \times S^{2n-2}$, whose boundary is diffeomorphic to $S^1\times S^{2n-2}$. Thus, the argument in the proof of Proposition \ref{proposition-index-localization} also  works in this setting and  we  obtain
   $$Ind(D^{\otimes L}_+)=(-1)^n\sum_{1\leq i\leq l}\deg^*(s_{x_i}).$$
\end{pf}

 \textbf{Acknowledgement}: First, the author  wants to show the best gratitude to  Professor  Huitao Feng, who brought this problem to the  focus. The author also would like to express the special thanks  to  Professor  Mikio Furuta, for the long time discussion and helping, and comments that greatly improved the manuscript. This research is supported by the Todai Fellowship.

\end{document}